\documentclass[a4paper,11pt]{article}

\usepackage{amsmath,amsfonts,amssymb,amsthm} 
\usepackage[explicit]{titlesec}
\usepackage{graphicx}
\usepackage{enumitem}
\setlist{nolistsep}
\usepackage{stmaryrd}
\usepackage{mathpple}
\usepackage{enumerate}
\usepackage{textcomp}
\usepackage[utf8]{inputenc}
\usepackage{mathtools}
\usepackage{epsfig}
\usepackage{xcolor}
\usepackage{array}
\usepackage{calc}
\usepackage[colorlinks=false]{hyperref}
\usepackage[noabbrev]{cleveref}
\usepackage{dsfont}
\usepackage{wrapfig}
\usepackage{fancyhdr}
\usepackage{float}
\usepackage[font={small}]{caption}
\usepackage{eurosym}
\usepackage{csquotes}
\usepackage{bm}
\usepackage[a4paper,left=3cm,right=3cm,top=2.5cm,bottom=2.5cm,footskip=1.8cm,headheight=0.5cm]{geometry}

\DeclareMathAlphabet{\mathcal}{OMS}{cmsy}{m}{n}
\newcommand{\dR}{\mathbb{R}}
\newcommand{\dP}{\mathbb{P}}

\newcommand{\E}{\mathbb{E}}

\newcommand{\cA}{\mathcal{A}}

\newcommand{\cH}{\mathcal{H}}

\newcommand{\cF}{\mathcal{F}}

\newcommand{\cI}{\mathcal{I}}
\newcommand{\cJ}{\mathcal{J}}
\newcommand{\cK}{\mathcal{K}}
\newcommand{\cR}{\mathcal{R}}
\newcommand{\cN}{\mathcal{N}}
\newcommand{\cM}{\mathcal{M}}

\newcommand{\cL}{\mathcal{L}}

\newcommand{\1}{\mathds{1}}

\newcommand{\drm}{\mathrm{d}}

\newcommand{\eps}{\varepsilon}

\newcommand{\wh}[1]{\widehat{#1}}

\font\dbdb=cmr10 scaled\magstep1
\font\calcal=cmsy10 scaled\magstep1
\def\build#1_#2^#3{\mathrel{\mathop{\kern 0pt#1}\limits_{#2}^{#3}}}
\def\liml{\build{\longrightarrow}_{}^{{\mbox{\calcal L}}}}
\def\limp{\build{\longrightarrow}_{}^{{\mbox{\dbdb P}}}}

\numberwithin{equation}{section}

\setitemize[1]{label=--,partopsep=\parskip,topsep=-\parskip,itemsep=0pt,parsep=0pt}
\titleformat{\chapter}{\fontfamily{phv}\selectfont\LARGE\bfseries}{\thesection}{1em}{\fontfamily{phv}\selectfont\LARGE\bfseries #1}
\titleformat{\section}{\centering\fontfamily{phv}\selectfont\bfseries}{\thesection}{1em}{\centering\fontfamily{phv}\selectfont\bfseries #1}
\titleformat{\subsection}{\fontfamily{phv}\selectfont\small\bfseries\centering}{\thesubsection}{1em}{\fontfamily{phv}\selectfont\small\bfseries #1}
\titleformat{\subsubsection}{\fontfamily{phv}\selectfont\footnotesize\bfseries\centering}{\thesubsubsection}{1em}{\fontfamily{phv}\footnotesize\selectfont\bfseries #1}

\renewcommand*\thesection{\arabic{section}}

\makeatletter
\renewcommand{\@seccntformat}[1]{\llap{{\csname the#1\endcsname}\hspace{1em}}}   

\newtheoremstyle{thmit}{10pt}{10pt}
{\normalfont\itshape}
{}
{\small\bf\fontfamily{phv}\selectfont}
{\;}
{0.25em}
{\small\fontfamily{phv}\selectfont\thmname{#1}\nobreakspace\footnotesize\thmnumber{#2}.
\thmnote{\nobreakspace\the\thm@notefont\fontfamily{phv}\selectfont\footnotesize\bfseries\--\nobreakspace#3}}

\newtheoremstyle{rmq}
{5pt}
{5pt}
{\normalfont}
{}
{\small\bf\fontfamily{phv}\selectfont}
{\;}
{0.25em}
{\small\fontfamily{phv}\selectfont\thmname{#1}\nobreakspace\footnotesize\thmnumber{\@ifnotempty{#1}{}\@upn{#2}}
\thmnote{\nobreakspace\the\thm@notefont\fontfamily{phv}\selectfont\footnotesize\bfseries--\nobreakspace#3.}}
\makeatother

\newcounter{count}
\numberwithin{count}{section}
\newcounter{alpha}

\theoremstyle{thmit}
\newtheorem{theorem}[count]{Theorem\small}

\newtheorem{remark}[count]{Remark\small}

\theoremstyle{rmq}

\renewenvironment{proof}[2]
{\paragraph{\fontfamily{phv}\selectfont\bfseries\small Proof of #1 \footnotesize #2.}}%
{\begin{flushright}
\qed
\end{flushright}}

\begin{document}
\pagestyle{fancy}

\begin{center}
\fontfamily{phv}\selectfont\bfseries
\Large
\uppercase{
How to estimate the memory of the Elephant Random Walk
\\[2ex]
\normalsize Bernard Bercu and 
Lucile Laulin
}
\end{center}

 \AtEndDocument{\bigskip{\footnotesize%
\textsc{Université de Bordeaux, Institut de Mathématiques de Bordeaux,
UMR 5251, 351 Cours de la Libération, 33405 Talence cedex, France.} \par  
\textit{E-mail adresses :} \href{mailto:bernard.bercu@math.u-bordeaux.fr}{\texttt{bernard.bercu@math.u-bordeaux.fr}} ;  
\href{mailto:lucile.laulin@math.u-bordeaux.fr}{\texttt{lucile.laulin@math.u-bordeaux.fr}}}} 

\vspace{2ex}
\centerline{
\begin{minipage}[c]{0.7\textwidth}
{\small\section*{Abstract}\vspace{-1ex}
We introduce an original way to estimate the memory parameter of the elephant random walk, 
a fascinating discrete time random walk on integers
having a complete memory of its entire history. Our estimator is nothing more than a quasi-maximum likelihood estimator, based on a second order Taylor 
approximation of the log-likelihood function. We show the almost sure convergence of our estimate in the diffusive, critical and superdiffusive regimes.
The local asymptotic normality of our statistical procedure is established in the diffusive regime, while the local asymptotic mixed normality
is proven in the superdiffusive regime. Asymptotic and exact confidence intervals as well as statistical tests are also provided. 
All our analysis relies on asymptotic results for martingales and the quadratic variations associated.
}\medskip\\
\small
{\bf MSC:} primary 60G50; secondary 60G42; 62M09\medskip
\\
{\bf Keywords : } Elephant random walk; Maximum likelihood estimation; Almost sure convergence; Local asymptotic normality
\end{minipage}
 }

\setlength{\parindent}{0pt}


\section{Introduction}\thispagestyle{empty}


The elephant random walk (ERW) is a fascinating discrete-time random walk on integers, which was introduced in the early 2000s
by Sch\"utz and Trimper \cite{Schutz2004}, in order to investigate how long-range memory affects the behavior of the random walk. 
It was referred to as the ERW in allusion to the traditional saying that elephants can always remember where they have been before. 
The ERW shows three different regimes depending on the location of its memory parameter and is defined as follows.
The elephant starts at the origin at time zero, $S_0 = 0$. For the first step, the elephant moves to the right at point $1$ with probability 
$q$ and to the left at point $-1$ with probability $1-q$ for some $q$ in $[0,1]$. The next steps are performed by choosing
at random an integer $k$ among the previous times $1,\ldots,n$. Then, the elephant moves exactly in 
the same direction as that of time $k$ with probability $p$ or to the oppositve direction with the probability $1-p$, where 
the parameter $p$ lies in $[0,1]$. In other words, for all $n \geq 1$, 
\begin{equation}
\label{STEPS}
   X_{n+1} = \left \{ \begin{array}{ccc}
    +X_{k} &\text{ with probability } & p, \vspace{2ex}\\
    -X_{k} &\text{ with probability } & 1-p.
   \end{array} \right.
\end{equation}
Therefore, the position of the ERW at time $n+1$ is given by
\begin{equation}
\label{POSERW}
S_{n+1}=S_{n}+X_{n+1}.
\end{equation}
The asymptotic behavior of the ERW is closely related to the value of $p$ called the memory of the ERW. 
The ERW is said to be diffusive if $0 \leq p < 3/4$, critical if $p=3/4$ and superdiffusive if $3/4< p \leq 1$. 
Whatever the value of $p$ in $[0,1]$, it has been shown that
\begin{equation}
\label{ASCVG}
 \lim_{n \rightarrow \infty} \frac{S_n}{n}=0 \hspace{1cm} \text{a.s.}
\end{equation}
Moreover, it has been proven in the diffusive regime $0 \leq p < 3/4$ that
\begin{equation}
\label{AND}
\frac{S_n}{\sqrt{n}} \liml \cN \Bigl(0, \frac{1}{3-4p}\Bigr),
\end{equation}
while in the critical regime $p= 3/4$ that
\begin{equation}
\label{ANC}
\frac{S_n}{\sqrt{n \log n}} \liml \cN (0, 1).
\end{equation}
We refer the reader to \cite{Baur2016}, \cite{Coletti2017}, \cite{ColettiN2017}, \cite{Bercu2018}
and to the recent contributions \cite{Bertoin2020}, \cite{Coletti2019}, \cite{Fan2020}, \cite{Takei2020}, \cite{Vazquez2019}. 
In the superdiffusive regime $3/4<p \leq 1$, it has been established by three different approaches \cite{Baur2016},
\cite{Bercu2018}, \cite{Coletti2017} that
\begin{equation}
\label{ASCVGS}
 \lim_{n \rightarrow \infty} \frac{S_n}{n^{2p-1}}=L \hspace{1cm} \text{a.s.}
\end{equation}
where $L$ is a non-degenerate random variable which is not Gaussian \cite{Bercu2018}. However, the fluctuation of the ERW around its limit $L$
are Gaussian \cite{Kubota2019} since, on the event $\{L^2 > 0 \}$,
\begin{equation}
\label{ANS}
\sqrt{n^{4p-3}} \Bigl(\frac{S_n}{n^{2p-1}} -L \Bigr) \liml \cN \Bigl(0, \frac{1}{4p-3}\Bigr).
\end{equation}
In this paper, we shall focus our attention on the parametric estimation the memory parameter $p$. To the best of our knowledge, no one has tackled 
this statistical analysis. Our estimator is explicitly given, for all $n \geq 2$, by 
\begin{equation}
\label{DEFPHAT}
 	\wh{p}_n  = \displaystyle \frac{\displaystyle\sum_{k=1}^{n-1} \frac{S_k}{k}\Bigl(X_{k+1} +\frac{S_k}{k}\Bigr)}{\displaystyle2\sum_{k=1}^{n-1}\Bigl(\frac{S_k}{k}\Bigr)^2}.
\end{equation}
The paper is organized as follows. In  Section \ref{S-QMLH}, we explain in detail how we are led to introduce the estimator $\wh{p}_n$
via a quasi-maximum likelihood approach. Section \ref{S-MR} is devoted to the main results of the paper.
We show the almost sure convergence of $\wh{p}_n$ to $p$ whatever the value of the memory parameter.
This preliminary estimation allows us to say whether the ERW is in the diffusive, critical or superdiffusive regimes.
The local asymptotic normality of our statistical procedure is established in the diffusive regime, while the local asymptotic mixed normality
is proven in the superdiffusive regime. In both regimes, asymptotic and exact confidence intervals as well as statistical tests are also provided.
Our martingale approach is described in Appendix A, while all technical proofs are postponed to Appendix B.


\section{Quasi-maximum likelihood estimation}
\label{S-QMLH}


Denote by $\cF_n = \sigma(X_1,\ldots,X_n)$ the $\sigma$-algebra of events occurring up to time $n$. It follows from
\eqref{STEPS} that for all $n \geq 1$, 
\begin{align*}
	\dP(X_{n+1} =1\mid\cF_n) & = \frac{p}{n} \sum_{k=1}^n\1_{\{X_k=1\}} + \frac{(1-p)}{n}\sum_{k=1}^n\1_{\{X_k=-1\}}, \\
	        & = \frac{p}{2n} \Bigl( n+S_n \Bigr) + \frac{(1-p)}{2n}\Bigl( n-S_n \Bigr), \\
			& = \frac{1}{2}\Big(1+(2p-1)\frac{S_n}{n}\Big).
\end{align*}
It clearly means that the conditional distribution of $X_{n+1}$ given $\cF_n$ is a Rademacher $\cR(p_n)$
distribution where
\begin{equation}
\label{DEFPNA}
p_n=\frac{1}{2}\Big(1+a\frac{S_n}{n}\Big)
\hspace{1cm}\text{and}\hspace{1cm} a=2p-1.
\end{equation}
Therefore, we obtain that for $x_{n+1}\in\{-1,1\}$
\begin{equation}
\label{Proba-Xn}
	\dP(X_{n+1} = x_{n+1}\mid\cF_n) = p_n^{(1+x_{n+1})/2}(1-p_n)^{(1-x_{n+1})/2}.
\end{equation}
For all $n \geq 1$ and $x \in \dR^n$ with $x=(x_1, \ldots,x_n)$, let 
$\dP_p(x)=\dP(X_1=x_1,\ldots,X_n=x_n)$. We clearly deduce from \eqref{Proba-Xn} that for all $n \geq 2$,
\begin{align*}
	\dP_{p}(x) & = \prod_{k=1}^{n-1} \dP(X_{k+1}=x_{k+1} \mid X_1=x_1, \ldots,X_{k}=x_{k}) \dP(X_1=x_1), \\
	 		& =  \prod_{k=1}^{n-1} \Big(p_k^{(1+x_{k+1})/2}(1-p_k)^{(1-x_{k+1})/2}\Big) q^{(1+x_{1})/2}(1-q)^{(1-x_{1})/2}
\end{align*}
where, for all $1 \leq k \leq n$, $S_k$ is replaced by $s_k=x_1+ \cdots+x_k$ in the definition of $p_k$. Consequently, the likelihood function associated
with $(X_1, \ldots, X_n)$ is given by
\begin{equation}
\label{LIKE}
	L_{n}(p) = \prod_{k=1}^{n-1} \Big(p_k^{(1+X_{k+1})/2}(1-p_k)^{(1-X_{k+1})/2}\Big) q^{(1+X_{1})/2}(1-q)^{(1-X_{1})/2}.
\end{equation}
It is easier to work with the log-likelihood function $\ell_{n}(p)= \log(L_{n}(p))$. We have from \eqref{LIKE} that
\begin{align}
	\ell_n(p)&  =   \sum_{k=1}^{n-1} \Big(\frac{1+X_{k+1}}{2}\Big)\log p_k + \Big(\frac{1-X_{k+1}}{2}\Big)\log (1-p_k) \nonumber \\
				&\quad + \Big(\frac{1+X_{1}}{2}\Big)\log q+ \Big(\frac{1-X_{1}}{2}\Big)\log (1-q).
\label{LOGLIKE}
\end{align}
Hence, if $\overline{X}_n$ stands for the empirical mean of $(X_1, \ldots, X_n)$, it follows
from \eqref{DEFPNA} and \eqref{LOGLIKE} that
\begin{align}
\label{EMVDEV}
	\ell_n'(p)&  =  \sum_{k=1}^{n-1} \big(1+X_{k+1}\big)\Big(\frac{\overline{X}_k}{1+a\overline{X}_k} \Big)- 
	\big(1-X_{k+1}\big) \Big(\frac{\overline{X}_k}{1-a\overline{X}_k} \Big),\nonumber \\
	& = \sum_{k=1}^{n-1} \frac{2\overline{X}_k(X_{k+1}-a\overline{X}_k)}{1-a^2\overline{X}_k^2},\nonumber\\
	& = \sum_{k=1}^{n-1} \frac{2X_{k+1}\overline{X}_k}{1+	aX_{k+1}\overline{X}_k}.
\end{align}
It is well-known that the process $(\ell_n'(p))$ is a locally square integrable martingale \cite{Heyde1975}. Its predictable quadratic variation is nothing else than
the conditional Fisher information $I_n(p)$ associated with $(X_1, \ldots, X_n)$. We shall see that
\begin{equation}
\label{CFISHER}
	I_{n}(p) = \sum_{k=1}^{n-1} \frac{\overline{X}_k^2}{p_k(1-p_k)} .
\end{equation}
It is not possible to find an explicit solution of the equation $\ell_n'(p)=0$. However, we already saw from \eqref{ASCVG} that
whatever the value of $p$ in $[0,1]$, $\overline{X}_n$ goes to zero almost surely. Consequently, it makes sense to replace $\ell_n(p)$ by
its second order Taylor approximation 
\begin{align}
	\lambda_n(p)	 &= \sum_{k=1}^{n-1} a\overline{X}_k \Big(X_{k+1} -\frac{a}{2}\overline{X}_k\Bigr) -(n-1) \log 2\nonumber \\
				& + \Big(\frac{1+X_{1}}{2}\Big)\log q+ \Big(\frac{1-X_{1}}{2}\Big)\log (1-q).
\label{EMVAPP}
\end{align}	
Since $a=2p-1$, \eqref{EMVAPP} clearly implies that
$$
\lambda_n'(p) =\sum_{k=1}^{n-1}2X_{k+1}\overline{X}_k \big(1 -aX_{k+1}\overline{X}_k\bigr)
\hspace{1cm}\text{and}\hspace{1cm} 
\lambda_n^{''}(p)=-4\sum_{k=1}^{n-1}\overline{X}_k^2.
$$
Therefore, $\lambda_n$ is a strictly concave function reaching its maximum at the value where its first derivative is equal to zero, which leads to 
\begin{equation*}
 	\wh{p}_n  = \displaystyle \frac{\displaystyle\sum_{k=1}^{n-1} \frac{S_k}{k}\Bigl(X_{k+1} +\frac{S_k}{k}\Bigr)}{\displaystyle2\sum_{k=1}^{n-1}\Bigl(\frac{S_k}{k}\Bigr)^2}.
\end{equation*}
It appears that our statistical approach is the most efficient strategy as it satisfies the local asymptotic normality (LAN) property in the diffusive regime
and the local asymptotic mixed normality (LAMN) property in the superdiffusive regime \cite{VdV1998}.


\section{Main results}
\label{S-MR}


Our first result deals with the almost sure convergence of $\wh{p}_n$ to $p$.

\begin{theorem}
\label{T-ASCVG}
Whatever the value of the memory parameter $p$ in $[0,1]$, $\wh{p}_n$ is a strongly consistent estimator of $p$,
\begin{equation}
	\label{ASCVGPN}
		\lim_{n\to\infty} \wh{p}_n = p \hspace{1.5cm} \text{a.s.}
	\end{equation}
\end{theorem}


\subsection{The diffusive regime}

Our next result is devoted to the asymptotic normality of the estimator $\wh{p}_n$ in the diffusive regime where $0 \leq p <3/4$.
Denote by $I(p)$ the asymptotic Fisher information
\begin{equation}
\label{AFISH-DR}
 I(p)= \frac{4}{3-4p}.
\end{equation}

\begin{theorem}
\label{T-AN-DR}
We have the asymptotic normality
\begin{equation}
\label{AN-DR}
	\sqrt{\log n}\big(\wh{p}_n -p\big)\liml \cN\Big(0,\frac{3-4p}{4}\Big).
\end{equation}
Its means that $\wh{p}_n$ is an asymptotically efficient estimator of $p$.
In particular,
\begin{equation}
\label{ANSLUT-DR}
	2\sqrt{\log n}\frac{\big(\wh{p}_n -p\big)}{\sqrt{3-4\wh{p}_n}}\liml \cN(0,1).
\end{equation}
\end{theorem}

We now focus our attention on the LAN property in the diffusive regime.

\begin{theorem}
\label{T-LAN-DR}
The sequence of experiments $(\dP_{n}(p), p\in [0,3/4[)$ is locally asymptotically normal. More precisely,
there exists a sequence of real random variables $(\Delta_n(p) )$ such that
$$
\Delta_n(p) \liml \cN\big(0,I(p)\big) 
$$
and for any sequence of real numbers $(h_n)$ converging to $h$, the log-likelihood ratio satisfies
\begin{equation}
\label{LAN-DR}
\log\Big(\frac{L_n(p+(\log n)^{-1/2}h_n)}{L_n(p)}\Big)=h\Delta_n(p) - \frac{h^2}{2}I(p) + o(1) \hspace{1cm} \text{a.s.}
\end{equation}
\end{theorem}

Our next result concerns an asymptotic confidence interval for the memory parameter $p$.
\begin{theorem} 
\label{T-CI-DR}
In the diffusive regime and for any $0<\alpha<1$, we have the asymptotic confidence interval for $p$ with confidence level $1-\alpha$,
\begin{equation}
\label{CI-DR}
	\cI(p)=\Bigl[\wh{p}_n-\frac{\sqrt{3-4\wh{p}_n}}{2\sqrt{\log n}}t_{1-\alpha/2},\ \wh{p}_n+\frac{\sqrt{3-4\wh{p}_n}}{2\sqrt{\log n}}t_{1-\alpha/2}\Bigr]
\end{equation}
where $t_{1-\alpha/2}$ stands for the $(1-\alpha/2)$-quantile of the standard $\cN(0,1)$ distribution.
\end{theorem}


\subsection{The critical regime}

We now focus our attention on the more complicated critical regime where $p=3/4$. 
Denote by $V_n$ a suitable approximation of the conditional Fisher information $I_n(p)$ given by \eqref{CFISHER},
\begin{equation}
\label{DEF-VN}
V_n= 4 \sum_{k=1}^{n-1}\Bigl(\frac{S_k}{k}\Bigr)^2.
\end{equation}

\begin{theorem}
\label{T-CR}
We have the  convergence in distribution
\begin{equation}
\label{VN-CR}
\frac{1}{(\log n)^2} V_n \liml 4\Lambda
\end{equation}
where $\Lambda$ stands for the integral of the squared standard Brownian motion
\begin{equation}
\label{DEF-LAMB}
\Lambda=\int_{0}^1 B_t^2 \drm t.
\end{equation}
\end{theorem}

\begin{remark}
It is impossible to prove the almost sure convergence as well as the convergence in probability in \eqref{VN-CR}. By the sharp analysis of Li \cite{Li1992PTRF, Li1992SPA} concerning the $L_2$-norm of the Brownian motion, we can only show that
$$
\liminf_{n\to\infty} \Bigl(\frac{\log \log \log n}{(\log n)^2} \Bigr)V_n = \frac{1}{2}
\hspace{1cm} \text{a.s.} 
$$
while
$$
\limsup_{n\to\infty} \Bigl( \frac{1}{ (\log n)^2 \log \log \log n} \Bigr) V_n = \frac{32}{\pi^2}
\hspace{1cm} \text{a.s.}
$$
This is the reason why we cannot establish the asymptotic normality of our estimator $\wh{p}_n$
in the critical regime.
\end{remark}

\begin{remark}
It follows from the Karhunen-Lo\`eve expansion of the Brownian motion that
\begin{equation}
\label{DECKL-LAMB}
\Lambda=\sum_{n=1}^\infty \frac{4}{(2n-1)^2 \pi^2} \xi_n^2
\end{equation}
where $(\xi_n)$ is a sequence of independent and identically distributed random variables with
$\cN(0,1)$ distribution, see e.g. Lemma 4 in \cite{Li1992PTRF}. Formula \eqref{DECKL-LAMB} allows the
numerical computation of the $\alpha$-quantiles of $\Lambda$, see \cite{Klein2003}.
\end{remark}


\subsection{The superdiffusive regime}

Our next result deals with the asymptotic normality of $\wh{p}_n$ 
in the superdiffusive regime where $3/4< p \leq 1$. 
We recall here that $L$ is the limiting non-degenerate random variable given in \eqref{ASCVGS}.

\begin{theorem}
\label{T-AN-SR}
Conditionally on the event $\{L^2>0\}$, we have the asymptotic normality
\begin{equation}
\label{AN-SR}
	\sqrt{V_n}(\wh{p}_n -p) \liml \cN(0,1).
\end{equation}
\end{theorem}
Figures 1 and 2 show the asymptotic normality of our estimator $\wh{p}_n$ in the diffusive and superdiffusive regimes with $p=0.4$ and $p=0.9$, respectively. The density function of the standard normal distribution is in red and the bins represent $N=3000$ different values of $\sqrt{V_n}(\wh{p}_n-p)$ for $n=1000$. We have used equation \eqref{AN-SR} to obtain both of the figures, as Theorem \ref{T-AN-SR} is also true in the diffusive regime. In fact, using directly the approximation of $V_n$ made in Theorem \ref{T-AN-DR} can not provide such good convergence results by simulations in the diffusive regime since $V_n$ increases almost surely to $4/(3-4p)$ with the slow speed $\log n$.
\begin{center}
\begin{minipage}[c]{.48\linewidth}
\begin{center}
\vspace{2ex} 
\label{Fig-TLC-04}
	\includegraphics[scale=0.44]{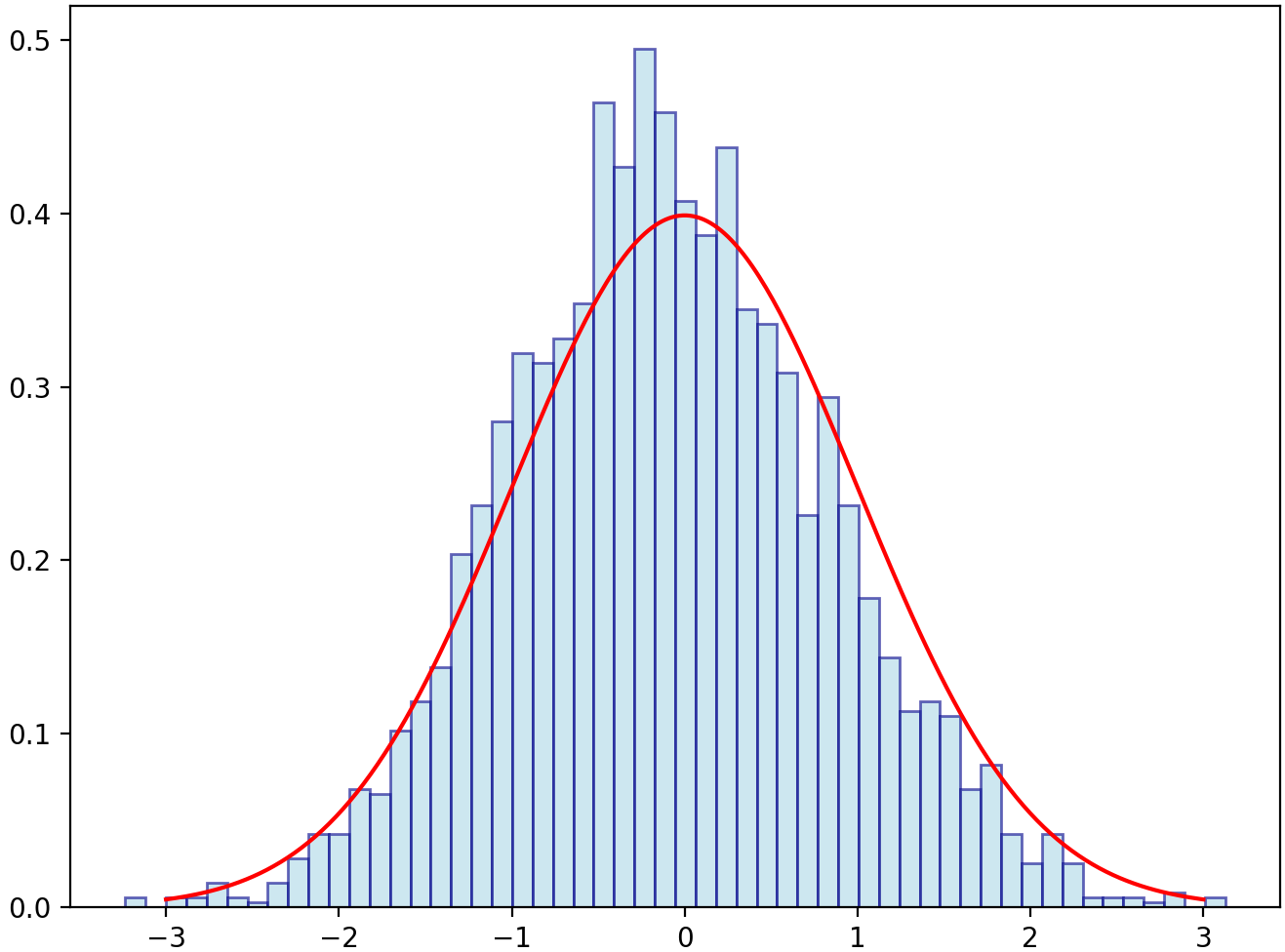}
\captionof{figure}{Asymptotic normality for $p=0.4$}
\end{center}
\end{minipage} \hfill %
\begin{minipage}[c]{.48\linewidth}
\begin{center}
\vspace{2ex} 
\label{Fig-TLC-09}
	\includegraphics[scale=0.44]{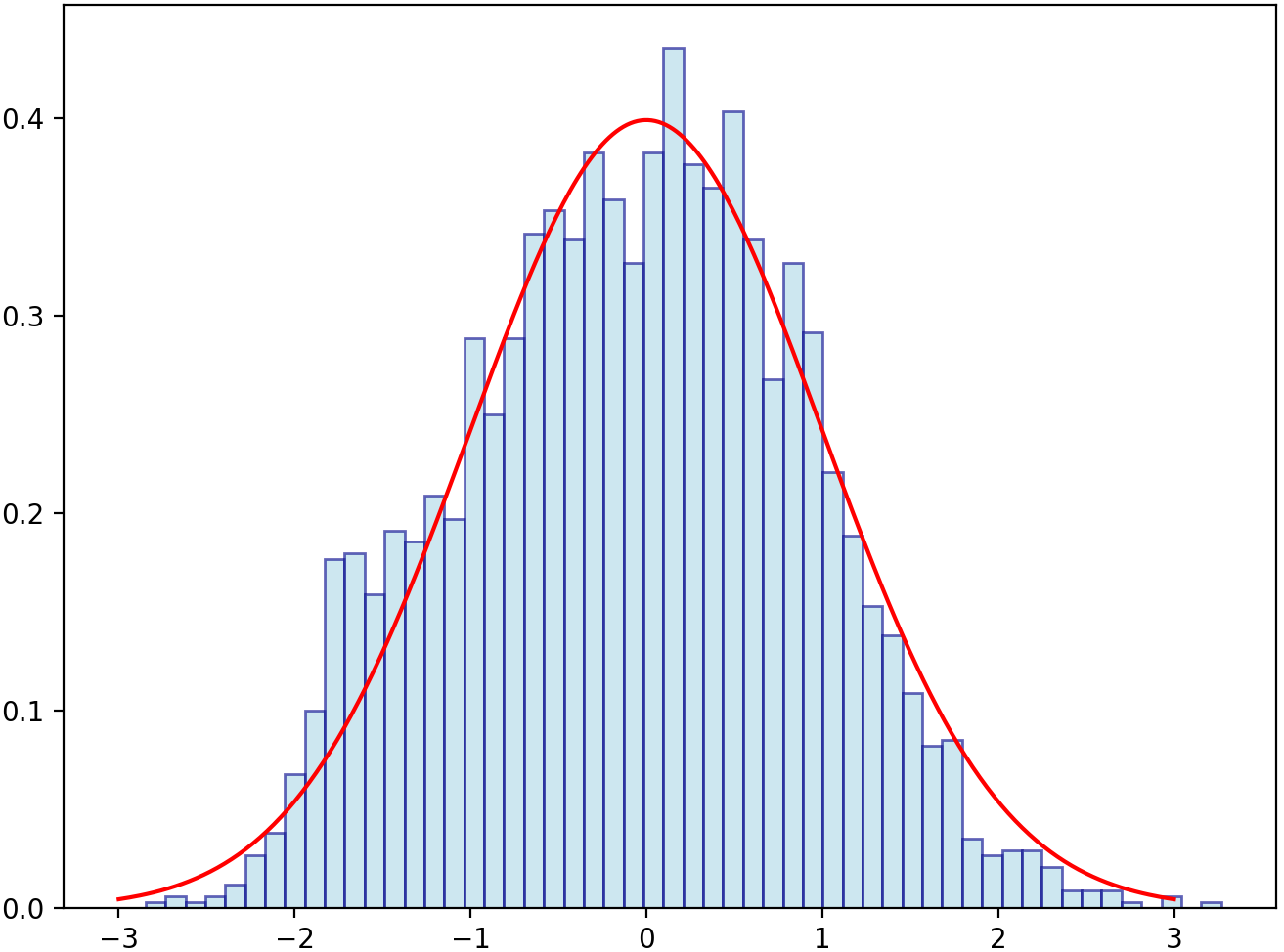}
\captionof{figure}{Asymptotic normality for $p=0.9$}
\end{center}
\end{minipage}	
\end{center}
\bigskip

The LAMN property in the superdiffusive regime is as follows.

\begin{theorem}
\label{T-LAMN-SR}
Conditionally on the event $\{L^2>0\}$, the sequence of experiments $(\dP_{n}(p), p\in ]3/4,1])$ is locally asymptotically mixed normal. More precisely,
there exists two sequences of real random variables $(\Delta_n(p) )$ and $(J_n(p))$ such that
\begin{equation*}
	\big(\Delta_n(p),J_n(p)\big) \liml \big(\Delta(p),J(p)\big)
\end{equation*}
and that the conditional distribution of $\Delta(p)$ given $J(p)=J$ is a standard $\cN(0,J)$ distribution,
and for any sequence of real numbers $(h_n)$ converging to $h$, the log-likelihood ratio satisfies
\begin{equation}
\label{LAMN-SR}
\log\Big(\frac{L_n(p+(n^{4p-3})^{-1/2}h_n)}{L_n(p)}\Big)=h\Delta_n(p) - \frac{h^2}{2}J_n(p) + o(1) \hspace{1cm} \text{a.s.}
\end{equation}
\end{theorem}

We also propose an asymptotic confidence interval for the memory parameter $p$.

\begin{theorem}
\label{T-CI-SR}
In the superdiffusive regime and for any $0<\alpha<1$, we have conditionally on the event $\{L^2>0\}$,
the asymptotic confidence interval for $p$ with confidence level $1-\alpha$,
\begin{equation}
\label{CI-SR}
	\cI(p)=\Bigl[\wh{p}_n-\frac{1}{\sqrt{V_n}}t_{1-\alpha/2},\ \wh{p}_n+\frac{1}{\sqrt{V_n}}t_{1-\alpha/2}\Bigr]
\end{equation}
where $t_{1-\alpha/2}$ stands for the $(1-\alpha/2)$-quantile of the standard $\cN(0,1)$ distribution.
\end{theorem}


\subsection{Exact confidence intervals}


Our purpose is now to provide an exact confidence interval for the  memory parameter $p$ whatever its value in $[0,1]$.
\begin{theorem}
\label{T-ECI}
For any $0<\alpha<1$, an exact confidence interval for $p$ with confidence level $1-\alpha$ is given, for all $n \geq 1$, by
\begin{equation} 
	\cJ(p)= \left[\wh{p}_n-\frac{2\sqrt{3n \log(2/\alpha)}}{V_n},\ \wh{p}_n+\frac{2\sqrt{3n \log(2/\alpha)}}{V_n}\right].
\end{equation}
Moreover, in the diffusive regime with $1/4 \leq p<3/4$, the exact confidence interval $\cJ(p)$ can be slightly improved by
\begin{equation}
		\cK(p)= \left[\wh{p}_n-\frac{\sqrt{29n \log(2/\alpha)}}{\sqrt{3}V_n},\ \wh{p}_n+\frac{\sqrt{29n \log(2/\alpha)}}{\sqrt{3}V_n}\right].
\end{equation}
\end{theorem}
\begin{remark}
Our confidence interval is better than the one obtained using Azuma-Hoeffding inequality which is given, for all $n \geq 3$, by
\begin{equation} 
	\cA(p)=\left[\wh{p}_n-\frac{2\sqrt{8n \log(2/\alpha)}}{V_n},\ \wh{p}_n+\frac{2\sqrt{8n \log(2/\alpha)}}{V_n}\right].
\end{equation}
\end{remark}

Figure 3 shows the three confidence intervals $\cI(p)$, $\cJ(p)$ and $\cA(p)$ in the superdiffusive regime
with $p=0.9$, for $n$ varying from $1$ to $100$. As expected, the asymptotic confidence interval $\cI(p)$ is
always more accurate than $\cJ(p)$ and $\cA(p)$, providing that the Gaussian approximation is justified.
One can also observe that $\cJ(p)$ and $\cA(p)$ are always true whatever the value of $n$ and that $\cJ(p)$
is more accurate than $\cA(p)$.

\begin{center}
\vspace{2ex} 
\label{Fig-IC}
	\includegraphics[scale=0.76]{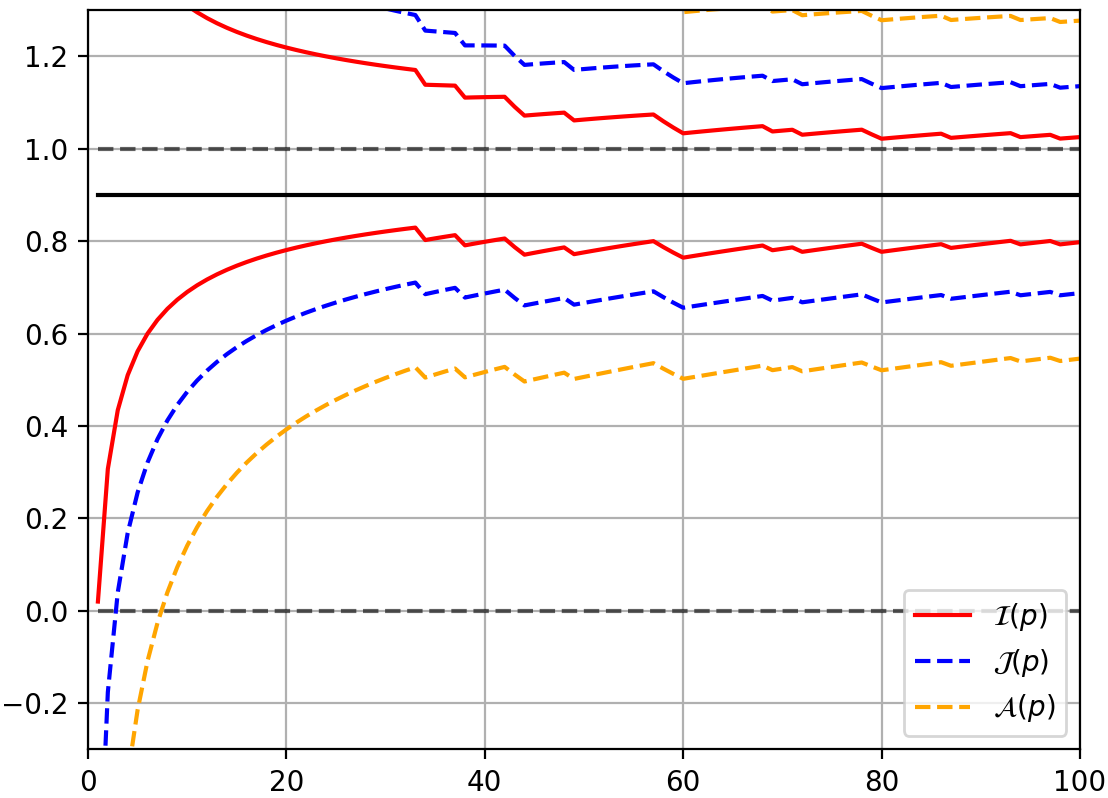}
\captionof{figure}{Confidence intervals for $p=0.9$ and $\alpha=0.05$}
\end{center}


\subsection{Statistical tests}

We are now in position to propose a bilateral statistical test built on our statistic $\wh{p}_n$. We start by fixing some memory value $0<p_0<1$ 
such that $p_0\neq 3/4$. Our goal is to test
\begin{equation*}
	\cH_0\ : \ ``p = p_0" \hspace{1cm}\text{against}\hspace{1cm} \cH_1\ : \ ``p \neq p_0".
\end{equation*}
\begin{theorem}
\label{T-TEST-p0}
Under the null hypothesis $\cH_0\ : \ ``p = p_0"$,
\begin{equation}
\label{CHI2}
	V_n(\wh{p}_n -p_0)^2 \liml \chi^2
\end{equation}
where $\chi^2$ has a Chi-square distribution with one degree of freedom.  
Moreover, under the alternative hypothesis $\cH_1\ : \ ``p \neq p_0"$,
\begin{equation}
	\lim_{n\to\infty} V_n(\wh{p}_n -p_0)^2 = +\infty \hspace{1cm} \text{a.s.}
\end{equation}
\end{theorem}
For a significance level $\alpha$ where $0<\alpha<1$, the acceptance and rejection regions are given by $\cA = [0,z_\alpha]$ and $\cR = ]z_\alpha,+\infty[$ where $z_\alpha$ stands for the $(1-\alpha)$-quantile of the Chi-square distribution with one degree of freedom.  The null hypothesis $\cH_0$ will not be rejected if the empirical value 
\begin{equation*}
	V_n(\wh{p}_n -p_0)^2 \leq z_\alpha
\end{equation*}
and will be rejected otherwise.
\vspace{1ex} \\
The purpose of our second test is to find out if the ERW is in the critical or the diffusive regime.
Concretely, we wish to test 
\begin{equation*}
	\cH_0\ : \ ``p = 3/4" \quad\text{against}\quad \cH_1\ : \ ``p < 3/4".
\end{equation*}
We immediatly obtain Theorem \ref{T-TEST-1}, whose proof directly follows from \eqref{VN-CR}.

\begin{theorem}
\label{T-TEST-1}
Under the null hypothesis $\cH_0\ : \ ``p = 3/4 "$,
\begin{equation}
	 \frac{1}{(\log n)^{2}}V_n \liml 4\Lambda 
\end{equation}
where $\Lambda$ is the integral of the squared Brownian motion given by
\eqref{DEF-LAMB}.
Moreover, under the alternative hypothesis $\cH_1\ : \ ``p < 3/4"$,
\begin{equation}
	\lim_{n\to\infty} \frac{1}{(\log n)^2}V_n = 0 \hspace{1cm} \text{a.s.}
\end{equation}
\end{theorem}

For a significance level $\alpha$ where $0<\alpha<1$, the acceptance and rejection regions are given by $\cA = [\lambda_{1-\alpha},+\infty[$ and $\cR = [0,\lambda_{1-\alpha}[$ where $\lambda_{1-\alpha}$ stands for the $\alpha$-quantile of 
the random variable $\Lambda$ which can be found in \cite{Klein2003}. For example, $\lambda_{0.05}=1.656$ and 
$\lambda_{0.10}=1.196$. The null hypothesis $\cH_0$ will not be rejected if the empirical value 
\begin{equation*}
	\frac{1}{4(\log n)^2} V_n \geq \lambda_{1-\alpha}
\end{equation*}
and will be rejected otherwise. The goal of our third is to find out
if the ERW is in the critical or superdiffusive regime. More precisely, we wish to test 
\begin{equation*}
	\cH_0\ : \ ``p = 3/4" \quad\text{against}\quad \cH_1\ : \ ``p > 3/4".
\end{equation*}

\begin{theorem}
\label{T-TEST-2}
Under the null hypothesis $\cH_0\ : \ ``p = 3/4 "$,
\begin{equation}
	 \frac{1}{(\log n)^{2}}V_n \liml 4\Lambda
\end{equation}
where $\Lambda$ is given by \eqref{DEF-LAMB}.
Moreover, under the alternative hypothesis $\cH_1\ : \ ``p > 3/4"$ and conditionnaly to $\{L^2>0\}$,
\begin{equation}
	\lim_{n\to\infty} \frac{1}{(\log n)^2}V_n = +\infty \hspace{1cm} \text{a.s.}
\end{equation}
\end{theorem}

For a significance level $\alpha$ where $0<\alpha<1$, the acceptance and rejection regions are given by $\cA = [0,\lambda_\alpha]$ and $\cR = ]\lambda_\alpha,+\infty[$ where $\lambda_\alpha$ stands for the $(1-\alpha)$-quantile of $\Lambda$, 
see \cite{Klein2003}. The null hypothesis $\cH_0$ will not be rejected if the empirical value 
\begin{equation*}
	\frac{1}{4(\log n)^2} V_n \leq \lambda_\alpha
\end{equation*}
and will be rejected otherwise. 
      

\renewcommand{\thesection}{\Alph{section}}
\renewcommand{\theequation}{A.\arabic{equation}}
\renewcommand{\thesubsection}{A.\arabic{subsection}}
\renewcommand{\thecount}{A.\arabic{count}}

\section*{Appendix A. Our martingal approach}
\addcontentsline{toc}{section}{Appendix A. The martingal approach}

\setcounter{section}{0}
\setcounter{equation}{0}
\setcounter{count}{0}


We already saw at the beginning of Section \ref{S-QMLH} that for all $n \geq 1$,
\begin{equation}
\label{CESTEPS}
 \E[X_{n+1}\mid\cF_n] = a \Big(\frac{S_n}{n}\Big)
\end{equation}
where $a=2p-1$. For all $n \geq 1$, let
\begin{equation*}
\eps_{n+1}= X_{n+1}-a\Big(\frac{S_n}{n}\Big)
\end{equation*}
with the initial value $\eps_1=X_1$. Since $(X_n)$ is a binary sequence of random variables taking values in
$\{+1,-1\}$, it clearly follows from \eqref{CESTEPS} that $(\eps_n)$ is a martingale difference sequence such that
for all $n \geq 1$, 
\begin{equation}
\label{CE2EPS}
 \E[\eps_{n+1}^2 \mid\cF_n] = 1- a^2 \Big(\frac{S_n}{n}\Big)^2.
\end{equation}
Equation \eqref{CE2EPS} immediately implies that
\begin{equation*}
\sup_{n \geq 1} \E[\eps_{n}^2] \leq 1.
\end{equation*}
Denote for all $n \geq 2$,
\begin{equation}
\label{DEF-MN}
	M_n = \sum_{k=1}^{n-1} \Big(\frac{S_k}{k}\Big)\eps_{k+1}
\end{equation} 
with $M_1=0$. As $|S_n| \leq n$, $(M_n)$ is a locally square integrable martingale. Its predictable quadratic variation is given by 
$\langle M\rangle_1 = 0$ and for all $n\geq 2$,
\begin{equation*}
 	\langle M\rangle_n = \sum_{k=1}^{n-1} \E[\Delta M_{k+1}^2 \mid\cF_{k}] = \sum_{k=1}^{n-1} \Big(\frac{S_k}{k}\Big)^2 \E[\eps_{k+1}^2\mid\cF_k].
 \end{equation*}
We obtain from \eqref{CE2EPS} that 
\begin{equation}
 	\langle M\rangle_n =  \sum_{k=1}^{n-1} \Big(\frac{S_k}{k}\Big)^2 \Big( 1- a^2 \Big(\frac{S_k}{k}\Big)^2\Big)= \sum_{k=1}^{n-1} \Big(\frac{S_k}{k}\Big)^2 - a^2 \sum_{k=1}^{n-1} \Big(\frac{S_k}{k}\Big)^4.
\label{QVMN}
 \end{equation}
Consequently, we deduce from \eqref{ASCVG} and \eqref{DEF-VN} that
 \begin{equation}
\label{ASCVGCMN}
\lim_{n \rightarrow \infty} \frac{V_n}{\langle M\rangle_n}=4 \hspace{1cm} \text{a.s.}
\end{equation} 
which means that the asymptotic behavior of the martingale $(M_n)$ is closely related to the one of the conditional Fisher information $I_n(p)$ and its approximation $V_n$.
Moreover, we have from \eqref{DEFPHAT} that
\begin{equation*}
 \wh{p}_n = \frac{\displaystyle\sum_{k=1}^{n-1} \frac{S_k}{k}\Big(X_{k+1}+\frac{S_k}{k}\Big)}{2\displaystyle\sum_{k=1}^{n-1} \Big(\frac{S_k}{k}\Big)^2} = 
 \frac{\displaystyle\sum_{k=1}^{n-1} \frac{S_k}{k}\Big(X_{k+1}-a\frac{S_k}{k}\Big)+(a+1)\sum_{k=1}^{n-1} \Big(\frac{S_k}{k}\Big)^2}{2\displaystyle\sum_{k=1}^{n-1} \Big(\frac{S_k}{k}\Big)^2}
\end{equation*}
which reduces, via \eqref{DEF-MN}, to
\begin{equation}
\label{DECMARTPN}
 \widehat{p}_n -p = \frac{2M_n}{V_n}.
\end{equation}
It ensures that the study of the asymptotic behavior of $\wh{p}_n$ can be achieved through convergence results for the martingale $(M_n)$.


\section*{Appendix B. Proofs of the main results}
\addcontentsline{toc}{section}{Appendix B. Proofs of the main results}

\renewcommand{\thesection}{\Alph{section}}
\renewcommand{\theequation}{B.\arabic{equation}}
\renewcommand{\thesubsection}{B.\arabic{subsection}}

\setcounter{section}{0}
\setcounter{subsection}{0}
\setcounter{equation}{0}
\setcounter{count}{0}


\begin{proof}{Theorem}{\ref{T-ASCVG}}
In the diffusive regime $0 \leq p < 3/4$, we have from the quadratic strong law given by Theorem 3.2 in \cite{Bercu2018} that
\begin{equation}
\label{QSL-DR}
\lim_{n \rightarrow \infty} \frac{1}{\log n} \sum_{k=1}^{n} \Big(\frac{S_k}{k}\Big)^2= \frac{1}{3-4p}
\hspace{1cm}\text{a.s.}
\end{equation}
which implies that
\begin{equation}
\label{VN-DR}
\lim_{n \rightarrow \infty} \frac{V_n}{\log n} = \frac{4}{3-4p}
\hspace{1cm}\text{a.s.}
\end{equation}
In the critical regime $p=3/4$, it follows once again from the quadratic strong law given by Theorem 3.5 in \cite{Bercu2018} that
\begin{equation*}
\lim_{n \rightarrow \infty} \frac{1}{\log \log n} \sum_{k=2}^{n} \Big(\frac{S_k}{k \log k}\Big)^2= 1
\hspace{1cm}\text{a.s.}
\end{equation*}
leading to
\begin{equation*}
\lim_{n \rightarrow \infty} V_n = + \infty
\hspace{1cm}\text{a.s.}
\end{equation*}

In the superdiffusive regime  $3/4< p \leq 1$, we deduce from \eqref{ASCVG} together with Toeplitz's lemma that
\begin{equation}
\label{QSL-SR}
\lim_{n \rightarrow \infty} \frac{1}{n^{4p-3}} \sum_{k=1}^{n} \Big(\frac{S_k}{k}\Big)^2= \frac{L^2}{4p-3}
\hspace{1cm}\text{a.s.}
\end{equation}
which ensures that
\begin{equation}
\label{VN-SR}
\lim_{n \rightarrow \infty} \frac{V_n}{n^{4p-3}} = \frac{4L^2}{4p-3}
\hspace{1cm}\text{a.s.}
\end{equation}
where $L$ is a non-degenrate random variable. Consequently, whatever the value of the memory parameter $p$ in $[0,1]$, we obtain that
$V_n$ increasing to infinity almost surely. Hence, we get from \eqref{ASCVGCMN} that $\langle M \rangle_n$ also goes to infinty almost surely  
in the three regimes.
Therefore, we can conclude from the strong law of large numbers 
for martingales given e.g. by Theorem 1.3.15 in \cite{Duflo1997} that
\begin{equation}
\label{ASCVGRATIOPN}
\lim_{n\rightarrow\infty} \frac{M_n}{V_n}=0\hspace{1cm} \text{a.s.}
\end{equation}
Finally, \eqref{DECMARTPN} together with \eqref{ASCVGRATIOPN} immediately lead to \eqref{ASCVGPN}.
\end{proof}


\subsection{The diffusive regime}


\begin{proof}{Theorem}{\ref{T-AN-DR}}
In the diffusive regime $0 \leq p < 3/4$, we already saw from \eqref{QSL-DR} that
\begin{equation*}
\lim_{n \rightarrow \infty} \frac{\langle M \rangle_n}{\log n}= \frac{1}{3-4p}
\hspace{1cm}\text{a.s.}
\end{equation*}
Moreover, $(M_n)$ satisfies the conditional Lindeberg condition, that is for all $\varepsilon >0$, 
\begin{equation*}
	\frac{1}{\log n}\sum_{k=1}^{n-1} \E\Big[|\Delta M_{k+1}|^2 \1_{|\Delta M_{k+1}| > \eps \sqrt{\log n}}\mid\cF_{k} \Big] \limp 0
\end{equation*} 
where, for all $n \geq 1$,
$$
\Delta M_{n+1}=M_{n+1}-M_n= \Big(\frac{S_n}{n}\Big)\eps_{n+1}.
$$
As a matter of fact, as $|S_n| \leq n$, we clearly have $|\eps_{n+1}| \leq 2$ and $| \Delta M_{n+1}| \leq 2$. Hence,  
we obtain that for all $\varepsilon >0$, 
\begin{align*}
	\frac{1}{\log n}\sum_{k=1}^{n-1} \E\Big[(\Delta M_{k+1})^2 \1_{|\Delta M_{k+1}| > \eps \sqrt{\log n}}\mid\cF_k \Big] 
	 & \leq  \frac{1}{\eps^2 (\log n)^2}\sum_{k=1}^{n-1} \E\Bigl[(\Delta M_{k+1})^4 \mid\cF_k  \Big], \\
	& \leq \frac{4}{\eps^2 \ (\log n)^2}\sum_{k=1}^{n-1} \Big(\frac{S_k}{k}\Big)^2, \\
	& \leq \frac{V_n}{\eps^2 (\log n)^2}.
\end{align*}
Therefore, we clearly deduce from \eqref{VN-DR} that 
\begin{equation*}
\lim_{n \rightarrow \infty} 	\frac{1}{\log n}\sum_{k=1}^{n-1} \E\Big[|\Delta M_{k+1}|^2 \1_{|\Delta M_{k+1}| > \eps \sqrt{\log n}}\mid\cF_{k} \Big] =0
	\hspace{1cm}\text{a.s.}
\end{equation*} 
which means that the conditional Lindeberg condition is satisfied. Hence, we can conclude from Corollary 3.1 in \cite{Hall1980} that 
\begin{equation}
\label{CVL-Mn}
	\frac{M_n}{\sqrt{\langle M\rangle_n}} \liml \cN(0,1).
\end{equation}
Finally, we obtain from \eqref{ASCVGCMN} and \eqref{DECMARTPN} together with Slutsky's Lemma that
\begin{equation*}
	\sqrt{V_n}(\wh{p}_n-p) \liml \cN(0,1)
\end{equation*}	
leading via \eqref{VN-DR} to
\begin{equation*}
	\sqrt{\log n}\big(\wh{p}_n -p\big)\liml \cN\Big(0,\frac{3-4p}{4}\Big).
\end{equation*}
One can observe that the asymptotic variance is the inverse of the Fisher information given by \eqref{AFISH-DR}, which completes the proof of Theorem
\ref{T-AN-DR}.

\end{proof}


\begin{proof}{Theorem}{\ref{T-LAN-DR}}


As in the proof of Theorem 7.2 in \cite{VdV1998} devoted to the Taylor expansion of the log-likelihood ratio, let
$$
\log(1+x)=x - \frac{x^2}{2}+x^2R(x)
$$ 
where the function $R(x)$ tends to zero as $x$ goes to zero. For any sequence of real numbers $(h_n)$ converging to $h$, we have
from \eqref{LOGLIKE} that
\begin{align*}
\ell_n(p+(\log n)^{-1/2}h_n)-\ell_n(p)&=
\sum_{k=1}^{n-1} \Big(\frac{1+X_{k+1}}{2}\Big)
	\log \Big(1+ \frac{2(\log n)^{-1/2}h_n\overline{X}_k}{1+a\overline{X}_k}\Big) \\
	& \quad 
	+  \sum_{k=1}^{n-1}\Big(\frac{1-X_{k+1}}{2}\Big)
	\log \Big(1- \frac{2(\log n)^{-1/2}h_n\overline{X}_k}{1-a\overline{X}_k}\Big)
\end{align*} 
Consequently, we obtain the Taylor expansion
\begin{align*}
\ell_n(p+(\log n)^{-1/2}h_n)-\ell_n(p)&=
\sum_{k=1}^{n-1} (1+X_{k+1}) \Big(\frac{(\log n)^{-1/2}h_n\overline{X}_k}{1+a\overline{X}_k}
	- \frac{(\log n)^{-1}h_n^2 \overline{X}_k^2}{(1+a\overline{X}_k)^2} \Big) \\
	& \hspace{-1cm} +2\sum_{k=1}^{n-1}(1+X_{k+1})\frac{(\log n)^{-1}h_n^2 \overline{X}_k^2}{(1+a\overline{X}_k)^2}  
	R\Big(\frac{2(\log n)^{-1/2}h_n\overline{X}_k}{1+a\overline{X}_k}\Big)\\
	& \hspace{-1cm} 
	-\sum_{k=1}^{n-1} (1-X_{k+1}) \Big(\frac{(\log n)^{-1/2}h_n\overline{X}_k}{1-a\overline{X}_k}
	+ \frac{(\log n)^{-1}h_n^2 \overline{X}_k^2}{(1-a\overline{X}_k)^2} \Big) \\
	& \hspace{-1cm} +2\sum_{k=1}^{n-1}(1-X_{k+1})\frac{(\log n)^{-1}h_n^2 \overline{X}_k^2}{(1-a\overline{X}_k)^2}  
	R\Big(\frac{2(\log n)^{-1/2}h_n\overline{X}_k}{1-a\overline{X}_k}\Big).
\end{align*} 
From now on, we are going to make repeated use that $(X_n)$ is a binary sequence of random variables taking values in
$\{+1,-1\}$. We can split the log-likelihood ratio into three terms,
\begin{equation}
\label{SPLITLLR-DR}
\ell_n(p+(\log n)^{-1/2}h_n)-\ell_n(p)= \frac{2h_n}{\sqrt{ \log n}}P_n -\frac{2h_n^2}{\log n}Q_n +\frac{2h_n^2}{\log n}R_n
\end{equation}
where
\begin{align*}
P_n & =  \frac{1}{2} \sum_{k=1}^{n-1} \Bigl(\frac{(1+X_{k+1})\overline{X}_k}{1+a\overline{X}_k} - \frac{(1-X_{k+1})\overline{X}_k}{1-a\overline{X}_k} \Bigr)
		 = \frac{1}{2} \sum_{k=1}^{n-1} \frac{2(X_{k+1}-a\overline{X}_k)\overline{X}_k }{1-(a\overline{X}_k)^2} \\
	  & =  \sum_{k=1}^{n-1} \frac{ (X_{k+1}-a\overline{X}_k)\overline{X}_k}{X_{k+1}^2-(a\overline{X}_k)^2}
		 = \sum_{k=1}^{n-1} \frac{ \overline{X}_k}{X_{k+1}+a\overline{X}_k} 
		 = \sum_{k=1}^{n-1} \frac{ X_{k+1}\overline{X}_k}{1+aX_{k+1}\overline{X}_k},
\end{align*}
\begin{align*}
Q_n & =  \frac{1}{2} \sum_{k=1}^{n-1} \Bigl(\frac{(1+X_{k+1})\overline{X}_k^2}{(1+a\overline{X}_k)^2} \!+\! \frac{(1-X_{k+1})\overline{X}_k^2}{(1-a\overline{X}_k)^2} \Bigr)
		 = \frac{1}{2} \sum_{k=1}^{n-1} \frac{2(1-2a\overline{X}_k X_{k+1}+a^2\overline{X}_k^2)\overline{X}_k^2}{(1-(a\overline{X}_k)^2)^2}  \\
	  & =  \sum_{k=1}^{n-1} \frac{ (X_{k+1}^2-2a\overline{X}_k X_{k+1}+a^2\overline{X}_k^2)\overline{X}_k^2}{( X_{k+1}^2-(a\overline{X}_k)^2)^2}
		 = \sum_{k=1}^{n-1} \frac{ ( X_{k+1}-a\overline{X}_k)^2\overline{X}_k^2}{( X_{k+1}+a\overline{X}_k)^2( X_{k+1}-a\overline{X}_k)^2} \\
		& =  \sum_{k=1}^{n-1} \frac{ \overline{X}_k^2}{(X_{k+1}+a\overline{X}_k)^2} = \sum_{k=1}^{n-1} \frac{ \overline{X}_k^2}{(1+aX_{k+1}\overline{X}_k)^2},
\end{align*}
and 
\begin{equation*}
	R_n  =  \sum_{k=1}^{n-1}\Bigl( \frac{(1+X_{k+1})\overline{X}_k^2}{(1+a\overline{X}_k)^2}  
	R\Big(\frac{2(\log n)^{-1/2}h_n\overline{X}_k}{1+a\overline{X}_k}\Big) \!+\!
	 \frac{(1-X_{k+1})\overline{X}_k^2}{(1-a\overline{X}_k)^2}  
	R\Big(\frac{2(\log n)^{-1/2}h_n\overline{X}_k}{1-a\overline{X}_k}\Big)\Bigr) .
\end{equation*}
On the one hand, we have
\begin{align*}
P_n & =   \sum_{k=1}^{n-1} \frac{ X_{k+1}\overline{X}_k}{1+aX_{k+1}\overline{X}_k}   =   \sum_{k=1}^{n-1} \frac{X_{k+1}\overline{X}_k \bigl(1-aX_{k+1}\overline{X}_k\bigr)}{\bigl(1+aX_{k+1}\overline{X}_k\bigr)\bigl(1-aX_{k+1}\overline{X}_k\bigr)} \\
      & =   \sum_{k=1}^{n-1} \frac{ \overline{X}_k\bigl(X_{k+1}-a\overline{X}_k\bigr) }{1-a^2\overline{X}_k^2} 
       =  \sum_{k=1}^{n-1}\frac{ \overline{X}_k \,\eps_{k+1} }{1-a^2\overline{X}_k^2}.
\end{align*}
It clearly means that the sequence $(P_n)$ is a square integrable martingale. We obtain from \eqref{CE2EPS} that the predictable quadratic variation
associated with $(P_n)$ is given by 
\begin{equation*}
 	\langle P\rangle_n = \sum_{k=1}^{n-1} \frac{ \overline{X}_k^2}{\big(1-a^2\overline{X}_k^2\big)^2}\E[\eps_{k+1}^2\mid\cF_k]=
 	\sum_{k=1}^{n-1} \frac{ \overline{X}_k^2}{1-a^2\overline{X}_k^2}.
 \end{equation*}
Hence, we immediately deduce from \eqref{ASCVG} and \eqref{QSL-DR} that
\begin{equation*}
\lim_{n \rightarrow \infty} \frac{\langle P \rangle_n}{\log n}= \frac{1}{3-4p}
\hspace{1cm}\text{a.s.}
\end{equation*}
Moreover, as it was previously done for the martingale $(M_n)$, one can easily check that $(P_n)$ satisfies the conditional Lindeberg condition. 
Consequently, it follows from Corollary 3.1 in \cite{Hall1980} that
\begin{equation}
\label{ANPN}
	\frac{P_n}{\sqrt{\log n}} \liml \cN\Big(0,\frac{1}{3-4p}\Big).
\end{equation}
On the other hand, we also have from \eqref{ASCVG} and \eqref{DEF-VN} that
\begin{equation}
Q_n  =   \sum_{k=1}^{n-1} \frac{ \overline{X}_k^2}{(1+aX_{k+1}\overline{X}_k)^2}  =  \frac{1}{4} V_n + o(V_n) \hspace{1cm} \text{a.s.}
\label{DEVQN}
\end{equation}
In the same way,
\begin{equation}
| R_n | =  o\left( \sum_{k=1}^{n-1} \frac{ \overline{X}_k^2}{(1+a\overline{X}_k)^2}\right)+o\left( \sum_{k=1}^{n-1} \frac{ \overline{X}_k^2}{(1-a\overline{X}_k)^2}\right)   =   o(V_n) \hspace{1cm} \text{a.s.}
\label{DEVRN}
\end{equation}
Finally, we obtain from the conjunction of \eqref{SPLITLLR-DR}, \eqref{ANPN}, \eqref{DEVQN} and \eqref{DEVRN} that
\begin{equation}
\label{LANF-DR}
\ell_n(p+(\log n)^{-1/2}h_n)-\ell_n(p)= h_n \Delta_n(p)  -\frac{h_n^2}{2}\frac{V_n}{\log n} +o(1) \hspace{1cm} \text{a.s.}
\end{equation}
where
$$
\Delta_n(p)=\frac{2 P_n}{\sqrt{\log n}} \liml \cN\Big(0,\frac{4}{3-4p}\Big),
$$
which is exactly what we wanted to prove.
\end{proof}


\begin{proof}{Theorem}{\ref{T-CI-DR}}
The proof directly follows from Theorem \ref{T-AN-DR}. Indeed, we obtain from the asymptotic normality \eqref{ANSLUT-DR} that
for any $0<\alpha<1$,
\begin{equation*} 
	\lim_{n\to\infty} \dP\Big(2\sqrt{\log n}\frac{\big|\wh{p}_n -p\big|}{\sqrt{3-4\wh{p}_n}}\leq t_{1-\alpha/2}\Big)
	 = 1 - \alpha
\end{equation*}
where $t_{1-\alpha/2}$ stands for the $(1-\alpha/2)$-quantile of the standard $\cN(0,1)$ distribution.
Moreover, one can easily see that
\begin{equation*} 
 \dP\Big(\!2\sqrt{\log n}\frac{\big|\wh{p}_n -p\big|}{\sqrt{3-4\wh{p}_n}}\!\!\leq t_{1-\alpha/2}\!\Big)
	\! = \! \dP\Big(\!\wh{p}_n-\frac{\sqrt{3-4\wh{p}_n}}{2\sqrt{\log n}}t_{1-\alpha/2}\leq p \leq
	\wh{p}_n+\frac{\sqrt{3-4\wh{p}_n}}{2\sqrt{\log n}}t_{1-\alpha/2}\!\Big).
\end{equation*}
It implies that 
\begin{equation*}
	\lim_{n\to\infty} \dP\Big(p \in\Big[\wh{p}_n-\frac{\sqrt{3-4\wh{p}_n}}{2\sqrt{\log n}}t_{1-\alpha/2},\ 
	\wh{p}_n+\frac{\sqrt{3-4\wh{p}_n}}{2\sqrt{\log n}}t_{1-\alpha/2}\Big]\Big) =1-\alpha,
\end{equation*}
which completes the proof of Theorem \ref{T-CI-DR}.
\end{proof}


\subsection{The critical regime}
\begin{proof}{Theorem}{\ref{T-CR}}
It follows from \eqref{POSERW} and \eqref{CESTEPS} with $a=1/2$ that for all $n \geq 1$,
\begin{equation*}
S_{n+1}=\Bigl(1+\frac{1}{2n}\Bigr) S_n + \varepsilon_{n+1}.
\end{equation*}
It clearly implies that for all $n \geq 1$,
\begin{equation}
\label{EMCR1}
\overline{X}_{n+1}=\Bigl(1-\frac{1}{2(n+1)}\Bigr) \overline{X}_n + \frac{1}{n+1}\varepsilon_{n+1}.
\end{equation}
Consequently, we obtain from \eqref{EMCR1} that for all $n \geq 2$,
\begin{equation*}
\overline{X}_{n}=\prod_{k=2}^n \Bigl(1-\frac{1}{2k}\Bigr) X_1 +
\sum_{k=2}^n \prod_{i=k+1}^n \Bigl(1-\frac{1}{2i}\Bigr) \frac{1}{k}\varepsilon_k
\end{equation*}
leading to
\begin{equation}
\overline{X}_{n}
= \frac{\Gamma(n+1/2)}{\Gamma(n+1)}\Bigl(  \frac{2}{\sqrt{\pi}} X_1 + \cM_n \Bigr)
\label{EMCR2}
\end{equation}
where
\begin{equation}
\cM_{n}=\sum_{k=2}^n \frac{\Gamma(k)}{\Gamma(k+1/2)}\varepsilon_k. 
\label{DEFCALM}
\end{equation}
We already saw that $(\eps_n)$ is a martingale difference sequence satisfying
\eqref{CE2EPS}. Hence, $(\cM_n)$ is a locally square integrable martingale with
predictable quadratic variation given, for all $n\geq 2$,
\begin{equation}
\label{CIPCALM}
 	\langle \cM\rangle_n = 
 	\sum_{k=2}^{n} \Big( \frac{\Gamma(k)}{\Gamma(k+1/2)} \Big)^{\!2}
 	-\frac{1}{4} \sum_{k=2}^{n} \Big( \frac{\Gamma(k)}{\Gamma(k+1/2)} \Big)^{\!2} \overline{X}_k^2.
 \end{equation}
Moreover, one can easily see that $\E[S_n^2]=n H_n$ where $H_n$ stands for the harmonic number
$$
H_n=\sum_{k=1}^n \frac{1}{k}.
$$
Therefore, we obtain from \eqref{CIPCALM} that
\begin{equation}
\label{CIPCALME}
 	\langle \cM\rangle_n - \E[\langle \cM\rangle_n] = 
 	-\frac{1}{4} \sum_{k=2}^{n} \Big( \frac{\Gamma(k)}{\Gamma(k+1/2)} \Big)^{\!2} 
 	\Big(\overline{X}_k^2- \frac{H_k}{k}\Big).
 \end{equation}
On the one hand, we have for all $n \geq 1$,
$$
n+\frac{1}{4}<\Big( \frac{\Gamma(n+1)}{\Gamma(n+1/2)} \Big)^{\!2}<n+\frac{1}{2},
$$
which is equivalent to
$$
\frac{1}{n}+\frac{1}{4n^2}<\Big( \frac{\Gamma(n)}{\Gamma(n+1/2)} \Big)^{\!2}<\frac{1}{n}+\frac{1}{2n^2}.
$$
It ensures that
\begin{equation}
\label{HCR}
\sum_{n=2}^{\infty} \Big( \frac{\Gamma(n)}{\Gamma(n+1/2)} \Big)^{\!2} \frac{H_n}{n} < 
\sum_{n=2}^{\infty} \frac{2 H_n}{n^2}<\infty.
\end{equation}
On the other hand, it follows from the quadratic strong law for the ERW given in Theorem 3.5 of
\cite{Bercu2018} that
\begin{equation}
\label{QSLCR}
\lim_{n\to\infty} \frac{1}{\log \log n} \sum_{k=2}^n \frac{\overline{X}_k^2}{(\log k)^2}=1 \hspace{1cm}\text{a.s.}
\end{equation}
Hence, we get from \eqref{QSLCR} together with Toeplitz lemme \cite{Duflo1997} that
\begin{equation*}
\lim_{n\to\infty} \frac{1}{\log \log n} \sum_{k=2}^n \frac{\overline{X}_k^2}{k}=0 \hspace{1cm}\text{a.s.}
\end{equation*}
which leads to
\begin{equation}
\label{QSLCRT}
\lim_{n\to\infty} \frac{1}{\log \log n} \sum_{k=2}^n \Big( \frac{\Gamma(k)}{\Gamma(k+1/2)} \Big)^{\!2}  \overline{X}_k^2=0 \hspace{1cm}\text{a.s.}
\end{equation}
Thus, we obtain from \eqref{CIPCALME}, \eqref{HCR} and \eqref{QSLCRT} that
$$
\langle \cM\rangle_n - \E[\langle \cM\rangle_n]=o(\log \log n) 
\hspace{1cm}\text{a.s.}
$$
Therefore, we deduce from the strong invariance principle for martingales given in
Theorem 2.1 of \cite{Shao1993} with $a_n=\log \log n$ and $b_n=\log n$ that
\begin{equation}
\label{SIPM}
\cM_n-B_{\log n}=o(\log \log n)  \hspace{1cm}\text{a.s.}
\end{equation}
Consequently, we obtain from \eqref{EMCR2} and \eqref{SIPM} the decomposition
\begin{equation}
\overline{X}_{n}
= \frac{\Gamma(n+1/2)}{\Gamma(n+1)} B_{\log n} +R_n \hspace{1cm}\text{a.s.}
\label{EMCR3}
\end{equation}
where the remainder $R_n$ 
satisfies
\begin{equation}
\sum_{k=1}^n R_{k}^2
= o\big(\log n (\log \log n)^2 \big)   \hspace{1cm}\text{a.s.}
\label{REMAINDERCR}
\end{equation}
In order to prove \eqref{VN-CR}, it only remains to show that
\begin{equation}
\label{CVGLAMBDA1}
\frac{1}{(\log n)^2} \sum_{k=1}^n \Big( \frac{\Gamma(k+1/2)}{\Gamma(k+1)}\Big)^2 B_{\log k}^2
\liml \Lambda
\end{equation}
where $\Lambda$ is the integral of the squared standard Brownian motion
\begin{equation*}
\Lambda=\int_{0}^1 B_t^2 \drm t.
\end{equation*}
We have for all $n \geq 1$,
$$
\frac{1}{n+1}<\Big( \frac{\Gamma(n+1/2)}{\Gamma(n+1)} \Big)^{\!2}<\frac{1}{n}.
$$
Consequently, the left-hand side in \eqref{CVGLAMBDA1} shares the same asymptotic behavior as
$$
\frac{1}{(\log n)^2} \sum_{k=1}^n  \frac{1}{k} B_{\log k}^2.
$$
Moreover, we have
\begin{align}
\sum_{k=1}^n  \frac{1}{k} B_{\log k}^2 &= \int_1^n \frac{1}{t}B^2_{\log t} \drm t + o(\log n) 
\hspace{1cm}\text{a.s.} \nonumber \\
&= \int_0^{\log n} B^2_{s} \drm s + o(\log n) 
\hspace{1cm}\text{a.s.}
\label{CVGLAMBDA2}
\end{align}
using the change of variables $s= \log t$. Hereafter, it follows from the self-similarity of the
Brownian motion that
\begin{equation}
\label{SELFS}
\int_0^{\log n} B^{2}_s \drm s \overset{\cL}{=}  (\log n)^{2} \int_{0}^{1} B_{t}^{2} \drm t
= (\log n)^2 \Lambda.
\end{equation}
Finally, we deduce from \eqref{EMCR3}, \eqref{REMAINDERCR}, \eqref{CVGLAMBDA2} and
\eqref{SELFS} that
$$
\frac{1}{(\log n)^2} \sum_{k=1}^n  \overline{X}_k^2 \liml \Lambda
$$
which completes the proof of Theorem \ref{T-CR}.

\end{proof}



\subsection{The superdiffusive regime}


\begin{proof}{Theorem}{\ref{T-AN-SR}}
In the superdiffusive regime $3/4<p\leq 1$, we already saw from \eqref{QSL-SR} that
\begin{equation*}
\lim_{n \rightarrow \infty} \frac{\langle M \rangle_n}{n^{4p-3}}= \frac{L^2}{4p-3}
\hspace{1cm}\text{a.s.}
\end{equation*}
\end{proof}
Moreover, as it was previously done in the diffusive regime, it is not hard to see that $(M_n)$ satisfies the conditional Lindeberg condition.
Hence, it follows from \eqref{DECMARTPN} together with Corollary 3.2 in \cite{Hall1980} that, conditionally on the event $\{L^2>0\}$, 
we have the asymptotic normality
\begin{equation*}
	\sqrt{V_n}(\wh{p}_n -p) \liml \cN(0,1).
\end{equation*}

\begin{proof}{Theorem}{\ref{T-LAMN-SR}}
As it was previously done in the proof of Theorem \ref{T-LAN-DR}, we can split the log-likelihood ratio into three terms,
\begin{equation}
\label{SPLITLLR-SD}
\ell_n(p+(n^{4p-3})^{-1/2}h_n)-\ell_n(p)= \frac{2h_n}{\sqrt{ n^{4p-3}}}P_n -\frac{2h_n^2}{n^{4p-3}}Q_n +\frac{2h_n^2}{n^{4p-3}}R_n \hspace{1cm}\text{a.s.}
\end{equation}
where the random variables $P_n$ and $Q_n$ are exactly the same, while the speed $\log n$ is replaced by $n^{4p-3}$ in the expression of $R_n$.
We immediately deduce from \eqref{ASCVGS} and \eqref{QSL-SR} that
\begin{equation*}
	\lim_{n\to\infty} \frac{\langle P \rangle_n}{n^{4p-3}} = \frac{L^2}{4p-3} \hspace{1cm}\text{a.s.}
\end{equation*}
Once again, $(P_n)$ satisfies the conditional Lindeberg condition in the superdiffusive regime. Hence, it follows from Corollary 3.2 in \cite{Hall1980} that, conditionally 
on the event $\{L^2>0\}$,
\begin{equation}
\label{ANPN-SD}
	\frac{P_n}{\sqrt{n^{4p-3}}}\liml L \times \cN\Big(0,\frac{1}{4p-3}\Big).
\end{equation}
Hereafter, we obtain from equations \eqref{SPLITLLR-SD}, \eqref{ANPN-SD}, \eqref{DEVQN} and \eqref{DEVRN}  that
\begin{equation}
\ell_n(p+(n^{4p-3})^{-1/2}h_n)-\ell_n(p)= h_n\Delta_n(p) -\frac{h_n^2}{2}J_n(p) +o(1) \hspace{1cm}\text{a.s.}
\end{equation}
where 
$$
\Delta_n(p) = \frac{2P_n}{\sqrt{ n^{4p-3}}} \hspace{1cm}\text{and} \hspace{1cm} J_n(p) = \frac{V_n}{n^{4p-3}}.
$$ 
Finally, as the convergence in \eqref{ANPN-SD} is stable \cite{Hall1980}, we immediatly have that
\begin{equation*}
 	\big(\Delta_n(p),J_n(p)\big) \liml \big(\Delta(p),J(p)\big)
 \end{equation*}
 where 
$$
J(p)= \frac{4 L^2}{4p-3}.
$$
In addition, conditionally on the event $\{J(p) = J\}$,
\begin{equation*}
	\Delta_n(p) = \frac{2P_n}{\sqrt{ n^{4p-3}}} \liml \cN\big(0,J\big),
\end{equation*}
which completes the proof of Theorem \ref{T-LAMN-SR}.
\end{proof}


\subsection{Exact confidence intervals}


\begin{proof}{Theorem}{\ref{T-ECI}}
In order to prove Theorem \ref{T-ECI}, we shall make use of concentration inequalities for martingales
\cite{Bercu2015}. 
First of all, one can observe that $(M_n)$ is a bounded difference martingale as equation \eqref{DEF-MN} implies
that for all $n \geq 2$,
\begin{equation*}
	\vert\Delta M_n\vert = \big \vert \eps_{n} \overline{X}_{n-1}\big\vert \leq (1+\vert a \vert)\big\vert\overline{X}_{n-1}\big\vert \hspace{1cm}\text{a.s.}
\end{equation*} 
Inspired by the Azuma-Hoeffding inequality for bounded difference martingales, denote
\begin{equation*}
	B_n = (1+\vert a \vert)^2 \sum_{k=1}^{n-1} \overline{X}_{k}^2.
\end{equation*}
Since $|a|\leq 1$, we clearly have from \eqref{QVMN} that
\begin{equation*}
5\langle M\rangle_n + B_n   \leq \frac{1}{4}\big(5+(1+\vert a \vert)^2\big)V_n  \leq 9n.
\end{equation*}
Hence, Theorem 3.4 in \cite{Bercu2015} ensures again that for any $x>0$,
\begin{equation}
\label{CIMN}
	\dP(|M_n|\geq  x)\leq 
			2 \exp \Big(-\frac{x^2}{3n}\Big).
\end{equation}
Consequently, it follows from \eqref{DECMARTPN} and \eqref{CIMN} that for any $x>0$,
\begin{equation}
\label{CIPN}
	\dP(V_n|\wh{p}_n-p|\geq 2n x)\leq 
			2 \exp \Big(-\frac{n x^2}{3}\Big).
\end{equation}
Hereafter, denote
$$
\alpha=2 \exp \Big(-\frac{n x^2}{3}\Big).
$$
As soon as $nx^2>3\log (2)$, the value $0<\alpha<1$. Therefore, we deduce from \eqref{CIPN} that
an exact confidence interval for $p$, with confidence level $1-\alpha$, is given by
\begin{equation*}
		\cJ(p)= \left[\wh{p}_n-\frac{2\sqrt{3n \log(2/\alpha)}}{V_n},\ \wh{p}_n+\frac{2\sqrt{3n \log(2/\alpha)}}{V_n}\right].
\end{equation*}
In the diffusive regime with $1/4 \leq p<3/4$, we have $|a| \leq 1/2$ which implies that
\begin{equation*}
5\langle M\rangle_n + B_n   \leq \frac{1}{4}\big(5+(1+\vert a \vert)^2\big)V_n  \leq  \frac{29}{4}n.
\end{equation*}
Hence, proceeding as in the previous calculation, we obtain  the exact confidence interval for $p$, with confidence level $1-\alpha$, 
\begin{equation*}
\cK(p)= \left[\wh{p}_n-\frac{\sqrt{29n \log(2/\alpha)}}{\sqrt{3}V_n},\ \wh{p}_n+\frac{\sqrt{29n \log(2/\alpha)}}{\sqrt{3}V_n}\right].
\end{equation*}
\end{proof}


\subsection{Statistical tests}

\begin{proof}{Theorem}{\ref{T-TEST-p0}}
The proof is quite straightforward. As a matter of fact, we already know from \eqref{AN-DR} or \eqref{AN-SR}
that under the null hypothesis $\cH_0$,
\begin{equation*}
	\sqrt{V_n}(\wh{p}_n -p_0) \liml \cN(0,1)
\end{equation*}
which immediatly implies \eqref{CHI2}.
It only remains to show that under the alternative hypothesis $\cH_1$, our test's statistic goes to infinity. Under $\cH_1$, we obtain from Theorem \ref{T-ASCVG}
that
	\begin{equation*}
	\lim_{n\to\infty} \wh{p}_n-p_0 = p-p_0 \hspace{1cm} \text{a.s.}
	\end{equation*}
and this limit is not zero. Consequently,
	\begin{equation*}
	\lim_{n\to\infty} V_n(\wh{p}_n-p_0)^2 = +\infty \hspace{1cm} \text{a.s.}
	\end{equation*}
as we already saw that whatever the value of the memory parameter $p$ in $[0,1]$, $V_n$ increasing to infinity almost surely,
completing the proof of Theorem \ref{T-TEST-p0}.
\end{proof}

\begin{proof}{Theorem}{\ref{T-TEST-1}}
The first part of Theorem \ref{T-TEST-1} immediatly follows from \eqref{VN-CR} under the null hypothesis.
It only remains to show that under the alternative hypothesis $\cH_1$, our test's statistic goes to 0. Under $\cH_1$, 
we have from equation \eqref{VN-DR}
that
	\begin{equation*}
	\lim_{n\to\infty} \frac{V_n}{\log n}=\frac{4}{3-4p} \hspace{1cm} \text{a.s.}
	\end{equation*}
which clearly implies that
\begin{equation*}
	\lim_{n\to\infty} \frac{V_n}{(\log n)^2}= 0 \hspace{1cm} \text{a.s.}
	\end{equation*}
which is exactly what we wanted to prove.
\end{proof}

\begin{proof}{Theorem}{\ref{T-TEST-2}}
The first part of Theorem \ref{T-TEST-2} follows once again from \eqref{VN-CR} under the null hypothesis.
It only remains to show that under the alternative hypothesis $\cH_1$, our test's statistic goes to infinity. 
Under $\cH_1$, we have from equation \eqref{VN-SR}
that
	\begin{equation*}
	\lim_{n\to\infty} \frac{V_n}{n^{4p-3}}=\frac{4L^2}{4p-3} \hspace{1cm} \text{a.s.}
	\end{equation*}
which implies that, conditionnaly on $\{L^2>0\}$,
\begin{equation*}
	\lim_{n\to\infty} \frac{V_n}{(\log n)^2}= +\infty \hspace{1cm} \text{a.s.}
	\end{equation*}
completing the proof Theorem \ref{T-TEST-2}.
\end{proof}


{\bf Acknowledgements. }The authors would like to thank Fabrice Gamboa and Florence Merlev\`ede for fruitful discussions on local asymptotic normality and strong invariance principles for martingales.
They also wish to thank Marc Arnaudon and Adrien Richou for sharing their in-depth
knowledge concerning functionals of Brownian motion.

\vspace{5ex}

{\small\bibliographystyle{acm}
\bibliography{bibliography}}

\end{document}